\documentclass{amsart}

\theoremstyle{definition}
\numberwithin{equation}{section}
\newcommand{\A}{\mathcal{A}}
\newcommand{\X}{\mathcal{X}}
\newcommand{\lra}{\longrightarrow}
\renewcommand{\theequation}{\arabic{equation}}
\begin{document}

\title{APPROXIMATE CONNES-AMENABILITY OF DUAL BANACH ALGEBRAS}

\author{G. H. Esslamzadeh}
\address{Department  of Mathematics, Faculty of Sciences, Shiraz University, Shiraz 71454, Iran, Telefax: (+98711)2281335}
\email{esslamz@shirazu.ac.ir}

\author{B. Shojaee}
\address{Department of Mathematics, Karaj Branch, Islamic Azad University, Karaj, Iran}
\email{shoujaei@kiau.ac.ir}

%\thanks{}

\begin{abstract}
We introduce the notions of approximate Connes-amenability and approximate strong Connes-amenability for dual Banach algebras. Then we characterize these two types of algebras in terms of approximate normal virtual diagonals and approximate $\sigma WC-$virtual diagonals. We investigate these properties for von Neumann algebras and measure algebras of locally compact groups. In particular we show that a von Neumann algebra is approximately Connes-amenable if and only if it has an approximate normal virtual diagonal. This is the ``approximate'' analog of the main result of Effros in [E. G. Effros, Amenability and virtual diagonals for von Neumann algebras, J. Funct. Anal. 78 (1988), 137-153].

We show that in general the concepts of approximate Connes-ameanbility and  Connes-ameanbility are distinct, but for measure algebras these two concepts coincide. Moreover cases where approximate Connes-amenability of $\A^{**}$ implies approximate Connes-amenability or approximate amenability of $\A$ are also discussed.

\end{abstract}

\maketitle

\baselineskip=18pt

\section{INTRODUCTION}

The concept of amenability for Banach algebras was introduced and studied for the first time by B. E. Johnson in [14]. Since then several variants of this concept have appeared in the literature each, as a kind of cohomological triviality. In [16], Johnson,  Kadison, and  Ringrose introduced a notion of amenability for von Neumann algebras which modified Johnson's original definition for Banach algebras in the sense that it takes the dual space structure of a von Neumann algebra into account. This notion of amenability was later called Connes-amenability by A. Ya. Helemskii [13]. Johnson in [15] showed that a Banach algebra $\A$ is amenable if and only if it has a virtual diagonal. A von Neumann algebraic analogue of this result was discovered by Haagerup [11]; See also [5], where the author introduces the notion of  normal virtual diagonal and presents another proof of Haagerup's result. Runde extended the notion of Connes-amenability to the larger class of dual Banach alge!
 bras [18] and studied certain concrete Banach algebras in the subsequent papers [20, 21, 22]. In particular he showed that existence of normal virtual diagonals implies Connes-amenability but the converse is no longer valid for arbitrary dual Banach algebras.

In all of the above mentioned concepts, all bounded derivations from a given Banach algebra $\A$ into certain Banach $\A$-bimodules  are required to be exactly inner. Gourdeau provided the following characterization of amenability; A Banach algebra $\A$ is amenable if and only if any bounded derivation from $\A$ into any Banach $\A$-bimodule is approximately inner, or equivalently weakly approximately inner [10, Proposition 2.1]. Motivated by Gourdeau's result, Ghahramani and Loy [8] introduced several approximate notions of amenability by requiring that all bounded derivations from a given Banach algebra $\A$ into certain Banach $\A$-bimodules to be approximately inner. However in contrast to Gourdeau's result, they removed the boundedness assumption on the net of implementing elements.  In the same paper and the subsequent one [9], the authors showed the distinction between each of these concepts and the corresponding classical notions and investigated properties of algebr!
 as in each of these new classes.
At the begining, Ghahramani and Loy asked which of the standard results on amenability work for the approximate concepts [See 8, page 233]; A question which identified the main direction of [6, 8, 9] and the present paper.

Motivated by the above question and [5], we introduce and study approximate Connes-amenability and approximate strong Connes-amenability. In Section 2 we present the definition and some basic properties of approximate Connes-amenability. An example presented at the begining of section 2, shows the distinction of Connes-amenability and approximate Connes-amenability. In Section 3  we introduce approximate strong Connes-amenability, approximate normal virtual diagonals and approximate $\sigma WC-$virtual diagonals. Then we show that a dual Banach algebra is approximately Connes-amenable [respectively, approximately strongly Connes-amenable] if and only if it has an approximate $\sigma WC-$virtual diagonal [respectively, approximately normal virtual diagonal]. In Section 4 which is the main part of this paper, we prove that a  von Neumann algebra is approximately Connes-amenable if and only if it has an approximate normal virtual diagonal. This is the ``approximate'' analog of !
 the main result of Effros [5]. In Section 5 we show that for a  locally compact group $G$, the measure algebra $M(G)$ of $G$ is  Connes-amenable if and only if it is approximately Connes-amenable if and only if it has an approximate normal virtual diagonal. This strengthens the main result of [21]. In the last section we show that under certain conditions approximate Connes-amenability of $\A^{**}$ implies approximate Connes-amenability or approximate amenability of $\A$.

We also should mention that some of our arguments were inspired by their classic analogs mostly from [5, 14, 18, 22].

Before proceeding further we recall some terminology.

Throughout $\A $ is  a Banach algebra and  $ \X $ is  a Banach $\A$-bimodule. Also the identity element of $\A$, whenever it exists, is denoted with $e$.  The dual space $\X^*$ of $\X$, is an $\A$-module, with  module actions
$$
\langle \phi .a\ ,\ x\rangle=\langle \phi\ ,\ a.x\rangle\quad ,\quad \langle a.\phi\ ,\ x\rangle=\langle \phi\ ,\ x.a\rangle\quad, \phi\in \X^*,\ x\in\X,\ a\in \A.
$$
 Using the natural $\A$-module structure of $\A^*$ the  first and second Arens multiplications on $\A^{**}$ that we denote by  ``.'' and ``$\Box$'' respectively, are defined by
$$
\langle m.n\ ,\ f\rangle =\langle m\ ,\ n.f\rangle\quad ,\quad\langle m\Box n\ ,\ f\rangle =\langle n\ ,\ f\Box m\rangle\quad f\in\A^*,\ m,n\in\A^{**}.
$$
The second dual of a Banach algebra, equipped with the first [respectively second] Arens product is a Banach algebra. We always consider the second dual of a Banach algebra with the first Arens product.

Throughout  ``derivation'' means ``bounded derivation'' and the set of all bounded derivations $D:\A \lra  \X $ is denoted by $Z^1(\A,\X)$. For $x\in  \X$ the map $ad_{x}(a)=a.x-x.a\quad (a\in \A )$ is called the inner derivation induced by $x$. A derivation $D:\A \lra  \X $ is approximately inner if there exists a net $(x_{\alpha})\subseteq  \X $ such that for every $a\in \A $, $D(a)=\lim_{\alpha}(a.x_{\alpha}-x_{\alpha}.a)$, the limit being in norm. We say that $\A $ is approximately amenable if for any $\A $-bimodule $ \X $, every derivation $D:\A \lra  \X ^{*}$ is approximately inner.

$\A $ is called a dual Banach algebra if there is a closed submodule $\A_{*}$ of $\A ^{*}$ such that $\A =(\A_{*})^{*}$.  In general the predual module is not necessarily unique. We will therefore assume that $\A$ always comes with a fixed predual $\A_*$. Measure algebras of  locally compact groups and second duals of Arens regular Banach algebras are examples of dual Banach algebras.

Let $\A $ be a dual Banach algebra and $\X$ be a Banach $\A $-bimodule. We call an element $\phi\in  \X ^{*}$ a normal element if
the  maps
$$
\A \lra  \X ^{*},\quad a\longmapsto\left\{\begin{array}{ll}a.\phi & \hbox{} \\
\phi.a & \hbox{}\end{array}\right.
$$
are $\omega^{*}-\omega^{*}$ continuous. If every element of $\X^*$ is normal,  then we say that $\X^*$ is normal. An element $x\in\X$ is called  $\omega^*$-weakly continuous if the module maps
$$
\A \lra  \X ,\quad a\longmapsto\left\{\begin{array}{ll}a.x & \hbox{} \\
x.a & \hbox{}\end{array}\right.
$$
are $\omega^*$-weakly continuous. The set of all $\omega^*$-weakly continuous elements of $\X$ is denoted by  $\sigma WC(\X)$.
A dual Banach algebra $\A$ is Connes-amenable if for every normal dual Banach $\A $-module $\X$, every $\omega^{*}-\omega^{*}$ continuous derivation $D\in Z^1(\A ,\X)$ is inner.

A left $\A$-submodule $\X$ of $\A^*$ is called left introverted if for every $\phi\in E$ and $m\in E^{*}$ the functional $m.\phi\in \A ^{*}$, which is defined by $\langle m.\phi,a\rangle=\langle m,\phi.a\rangle\quad (a\in \A )$, lies again in $E$. This turns $E^{*}$ into a dual Banach algebra by letting
$\langle nm,\phi\rangle=\langle n,m.\phi\rangle\quad(n,m\in E^{*},\phi\in E)$.
 The $WAP(\A ^{*})=\{\phi\in \A ^{*}: a\longmapsto a.\phi $ is weakly compact
  $, a\in \A \} $ is a left introverted subspace of $\A ^{*}$ and hence $WAP(\A ^{*})^{*}$ with the above product is a dual Banach algebra.

\section{DEFINITION AND BASIC PROPERTIES}

\textbf{Definition.} A dual Banach algebra $\A$ is approximately
Connes-amenable if for every normal, dual Banach $\A$-bimodule $\X$,
every $\omega^{*}-\omega^{*}$ continuous derivation $D\in Z^{1}(\A , \X )$ is
approximately inner.

The following example shows the distinction between  Connes-amenability and approximate Connes-amenability.

\textbf{Example 2.1.} Let $S$ be the set $\Bbb N$ of natural numbers with the binary operation $(m,n)\mapsto\max\{ m,n\}$. Then $S$ is a unital, commutative, weakly cancellative semigroup, that is, for every $s,t\in S$ the set $\{ x\in S\ :\ sx=t\}$ is  finite. Let $\A=\ell^1(S)$. Since $S$ is weakly cancellative, then by [3, Theorem 4.6] $\A$ is a dual Banach algebra with respect to the predual $c_0(S)$. If $\A$ is Connes-amenable, then by [4, Theorem 5.13] $S$ should be a group which is not the case. Thus $\A$ is not Connes-amenable. However as it was shown in [3, Example 10.10], $\A$ is approximately amenable and hence is approximately Connes-amenable.

\textbf{Proposition 2.2.} Suppose that $\A$ is approximately Connes-amenable. Then $\A$
has left and right approximate identities. In particular $\A ^{2}$ is
dense in $\A$.

{\it Proof.} Let $\X$ be the Banach $\A$-bimodule whose underlying linear space
is $\A$ equipped with the  module operations $a.x=ax$ and $x.a=0$, $(a\in \A,\ \ x\in  \X ).$

Obviously $\X$ is a normal dual Banach $\A$-bimodule and the identity map
on $\A$ is a $\omega^{*}-\omega^{*}$ continuous derivation. Since $\A$ is
approximately Connes-amenable, then there exists a net
$(a_{\alpha})\subseteq  \X $ such that $a=\lim_{\alpha}aa_{\alpha}\quad\quad(a\in \A).$

This means that $\A$ has a right approximate identity. Similarly, one
see that $\A$ has a left approximate identity.$\qed$

Let $(\A,\A_{*})$ be a dual Banach algebra, and let $\A^{\#}$ be the
Banach algebra $\A\oplus \Bbb C$. Then $\A^{\#}$ is a dual Banach algebra
with predual $\A_{*}\oplus \Bbb C$ and norm
$$\parallel(\mu,\alpha)\parallel=\max(\|\mu\|,|\alpha|)\quad\quad (\mu\in \A_{*},\alpha\in \Bbb C).$$

\textbf{Proposition 2.3.} Let $\A$ be a dual Banach algebra. $\A$ is
approximately Connes-amenable if and only if $\A^{\#}$ is
approximately Connes-amenable.

{\it Proof.} Let $D:\A^{\#}\lra  \X$ be a
$\omega^{*}-\omega^{*}$ continuous derivation where $\X$ is a normal dual Banach
$\A^{\#}$-bimodule. By [8, Lemma 2.3], $D=D_{1}+ad_{\eta}$ where
$D_{1}:\A^{\#}\lra e. \X .e$ is a $\omega^{*}-\omega^{*}$ continuous
derivation and $\eta\in  \X $. Since $e. \X .e$ is a normal dual
Banach $\A$-bimodule, then $D_{1}(e)=0$ and $D_{1}|_{\A}$ is approximately
inner; whence D is approximately inner. Thus $\A^{\#}$ is
approximately Connes-amenable.

Now suppose $D:\A\lra  \X$ is a $\omega^{*}-\omega^{*}$ continuous
derivation where  $ \X $ is a normal dual Banach $\A$-bimodule. Set
$$\widetilde{D}:\A^{\#}\lra  \X,\quad \widetilde{D}(a+\lambda e)=Da\quad (a\in \A, \lambda\in \Bbb C).$$

If we define $e.x=x.e=x\quad (e\in \A^{\#},x\in  \X )$, then $ \X $ turns into a normal dual Banach $\A^{\#}$-bimodule and
$\widetilde{D}$ is a $\omega^{*}-\omega^{*}$ continuous derivation. So
$\widetilde{D}$ is approximately inner, and hence so is D. It follows that
 $\A$ is approximately Connes-amenable.$\qed$

\textbf{Proposition 2.4.} Suppose that $\A$ is a dual Banach algebra with identity. Then $\A$ is approximately Connes-amenable if and only if every $\omega^{*}-\omega^{*}$ continuous derivation into every unital normal dual Banach $\A$-bimodule $\X$ is approximately inner.

{\it Proof.}  Suppose $D\in Z^{1}(\A, \X)$ is a
$\omega^{*}-\omega^*$ continuous derivation into the normal dual Banach bimodule
$\X$. By [8, Lemma 2.3], we have $D=D_{1}+ad_{\eta}$ where
$D_{1}:\A\lra e. \X .e$ is a derivation and $\eta\in  \X $.
Since D is a $\omega^{*}-$continuous derivation and $\X$  is a normal
dual Banach bimodule then $D_{1}$ is $\omega^{*}-$continuous and
$e. \X .e$ is normal. So by assumption $D_{1}$ is approximately inner,
and therefore $\A$ is approximately Connes-amenable. The converse holds obviously. $\qed$

\section{APPROXIMATE NORMAL VIRTUAL DIAGONALS}

Throughout this section we assume that $\A$ is a dual Banach algebra with identity. See Remark 3.4 at the end of this section regarding the non-unital case.

Let $L^2(\A,\Bbb C)$ be the space of all bounded bilinear functionals on $\A$ and $L_{\omega^{*}}^{2}(\A,\Bbb C)$ be the space of separately $\omega^{*}$ continuous elements of $L^2(\A,\Bbb C)$. Following the terminology of [5, 16], we turn $L^{2}(\A,\Bbb C)$ into a Banach $\A$-bimodule through the identification $L^{2}(\A,\Bbb C)\simeq (\A\widehat{\otimes}\A)^{*}$. Then the module actions of $\A$ on $L^{2}(\A,\Bbb C)$ are as follow.
$$
(a.F)(b,c)=F(b,ca),\quad (F.a)(b,c)=F(ab,c),\quad a,b,c\in\A ,\quad F\in L^2(\A ,\Bbb C).
$$
Clearly, $L_{\omega^{*}}^{2}(\A,\Bbb C)$ is a Banach $\A$-submodule of
$L^{2}(\A,\Bbb C)$.  Moreover we have a natural $\A$-bimodule map
$$
\theta:\A\otimes \A\lra L_{\omega^{*}}^{2}(\A,\Bbb C)^{*},\quad \theta(a\otimes b)(F)=F(a,b).
$$
 Since $\A_{*}\otimes\A_{*}\subseteq L_{\omega^{*}}^{2}(\A,\Bbb C)$ and $\A_{*}\otimes \A_{*}$
separates points of $\A\otimes \A$, then $\theta$ is one-to-one. We
will identify $\A\otimes \A$ with its image, writing
$$
\A\otimes \A\subseteq L_{\omega^{*}}^{2}(\A,\Bbb C)^{*}.
$$
The map $\Delta_{\A}$ is defined as follows.
$$
\Delta_{\A}:\A\widehat{\otimes} \A\lra \A,\quad\quad a\otimes b\longmapsto ab\quad (a,b\in \A).
$$
Since multiplication in a dual Banach algebra is separately
$\omega^{*}-\omega^{*}$-continuous, we have
$$
\Delta_{\A}^{*}(\A_{*})\subset L_{\omega^{*}}^{2}(\A,\Bbb C).
$$
So the restriction of $\Delta_{\A}^{**}$ to $L_{\omega^{*}}^{2}(\A,\Bbb C)^*$ turns into a Banach $\A$-bimodule homomorphism
$$
\Delta_{\omega^{*}}:L_{\omega^{*}}^{2}(\A,\Bbb C)^{*}\lra \A .
$$
Suppose $F\in L_{\omega^{*}}^{2}(\A,\Bbb C)$ and $M\in
L_{\omega^{*}}^{2}(\A,\Bbb C)^{*}$. We use the notation,
$$
\int F(a,b)dM(a,b)=\int FdM:=\langle M,F\rangle.
$$
More generally given a dual Banach space $ \X ^{*}$ and a bounded
bilinear function $F:\A\times \A\lra  \X ^{*}$ such that
$a\lra F(a,b)$ and $b\lra F(a,b)$ are
$\omega^{*}-\omega^{*}$-continuous, $\int FdM\in  \X ^{*}$ is defined  by
$$
\langle\int FdM,x\rangle=\int\langle F(a,b),x\rangle dM(a,b)\quad (x\in  \X ).
$$
Sometimes we also use the term $\int F(a,b)dM(a,b)$ for $\int FdM$.

\textbf{Definition.} A net $(M_{\alpha})$ in $L_{\omega^{*}}^{2}(\A,\Bbb C)^{*}$ is called an
approximate normal, virtual diagonal for $\A$ if for every $a\in\A$
$$
a.M_{\alpha}-M_{\alpha}.a\lra0\quad and \quad \Delta_{\omega^{*}}(M_{\alpha})\lra e,
$$
the limits being taken in norm.

It is well known that every dual Banach algebra with a normal virtual diagonal is Connes-amenable [18]. In the following theorem we extend this result to approximate Connes-amenability.

\textbf{Theorem 3.1.} If $\A$ has an approximate normal, virtual diagonal $\{M_{\alpha}\}$, then $\A$ is
approximately Connes-amenable.

{\it Proof.} Suppose $ \X $ is a normal dual Banach
$\A$-bimodule with predual $ \X _{*}$ and $D\in Z^{1}(\A, \X )$
 is $\omega^{*}-\omega^{*}$-continuous. Since $\A$ has an identity,
by Proposition 2.4 we can assume that $\X$  is unital. Since the bilinear map
$$
F:\A\times\A\lra\X,\quad F(a,b)=Da.b
$$
is separately $\omega^{*}-\omega^{*}$ continuous, then by the preceding remark we may define
$$
\phi_{\alpha}=\int F(a,b)dM_{\alpha}(a,b)=\int Da.b\, dM_{\alpha}\in  \X .
$$
For $c\in \A$, $x\in  \X _{*}$ we have
$$
\langle c.\phi_{\alpha},x\rangle=\langle\phi_{\alpha},x.c\rangle=\int\langle c.Da.b,x\rangle dM_{\alpha}(a,b)=\langle\int c.Da.b\, dM_{\alpha}(a,b),x\rangle.
$$
Therefore

\begin{eqnarray}
c.\phi_{\alpha}=\int c.Da.b\, dM_{\alpha}(a,b) \label{1}
\end{eqnarray}

and similarly

\begin{eqnarray}
\phi_{\alpha}.c=\int Da.bc\, dM_{\alpha}(a,b). \label {2}
\end{eqnarray}

So if we define $F_{x}\in L_{\omega^{*}}^{2}(\A,\Bbb C)$ by  $F_{x}(a,b)=\langle Da.b,x\rangle$, then the following relations hold.

\begin{eqnarray}
\int\langle D(ca).b,x\rangle dM_{\alpha}(a,b)=\int F_x.c(a,b) dM_{\alpha}(a,b)=\langle c.M_{\alpha},F_{x}\rangle \label{3}
\end{eqnarray}

\begin{eqnarray}
\int\langle Da.bc,x\rangle dM_{\alpha}(a,b)=\int c.F_x(a,b) dM_{\alpha}(a,b)=\langle M_{\alpha}.c,F_{x}\rangle\label{4}
\end{eqnarray}
By (3) and (4) we have
$$
|\langle\int D(ca).b\, dM_{\alpha}(a,b)-\int Da.bc\, dM_{\alpha}(a,b),x\rangle|\leq\|c.M_{\alpha}-M_{\alpha}.c\|\|F_{x}\| .
$$
So
\begin{eqnarray}
\|\int Dca.b\, dM_{\alpha}(a,b)-\int Da.bc\, dM_{\alpha}(a,b)\|\leq\|c.M_{\alpha}-M_{\alpha}.c\|\|D\|\|a\|\|b\|.\label{5}
\end{eqnarray}

If we define $G: \A\times\A\lra L_{\omega^{*}}^{2}(\A,\Bbb C)^*$ by $G(a,b)=a\otimes b$, then for every $F$ in $L_{\omega^{*}}^{2}(\A,\Bbb C)$ we have
$$
\langle \int GdM_{\alpha},F\rangle=\int\langle G(a,b),F\rangle dM_{\alpha}(a,b)
=\int F(a,b)dM_{\alpha}(a,b)=\langle M_{\alpha},F\rangle .
$$
So $M_{\alpha}=\int (a\otimes b)dM_{\alpha}(a,b)$. Now for every $t\in\A_*$,
$$
\langle \Delta_{\omega^{*}}(M_{\alpha}),t\rangle=\langle M_{\alpha},\Delta_{\A}^{*}(t)\rangle=
\int \langle a\otimes b ,\Delta_{\A}^{*}(t)\rangle dM_{\alpha}(a,b)=\langle\int ab\, dM_\alpha(a,b), t\rangle .
$$
Thus
\begin{eqnarray}
\Delta_{\omega^{*}}(M_{\alpha})=\int ab\, dM_{\alpha}(a,b). \label{6}
\end{eqnarray}
Moreover we have
$$
\aligned \langle Dc.\int ab\, dM_{\alpha}(a,b), x\rangle &=\int \langle ab, x.Dc\rangle dM_{\alpha}(a,b)\\
&=\int \langle Dc. ab, x\rangle dM_{\alpha}(a,b)=\langle \int Dc.ab\, dM_{\alpha}(a,b), x\rangle .\endaligned
$$
Therefore
\begin{eqnarray}
Dc.\int ab\, dM_{\alpha}(a,b)=\int Dc.ab\, dM_{\alpha}(a,b) . \label{7}
\end{eqnarray}
Now by (1), (2) and (7),
$$
\aligned c.\phi_{\alpha}-\phi_{\alpha}.c&=\int D(ca).b\, dM_{\alpha}(a,b)-\int Dc.ab\, dM_{\alpha}(a,b)
-\int Da.bc\, dM_{\alpha}(a,b)\\
&=\int D(ca).b\, dM_{\alpha}(a,b)-\int Da.bc\, dM_{\alpha}(a,b)- Dc.\int ab\, dM_{\alpha}(a,b). \endaligned
$$
Applying our assumption and (5) to the above identity  shows that
$$
lim_{\alpha}(\phi_{\alpha}.c-c.\phi_{\alpha})=Dc\quad (c\in \A).
$$
Therefore $D$ is approximately inner, and hence $\A$ is approximately
Connes-amenable. $\qed$

We don't know whether the converse of Theorem 3.1 is true in general. However we show in Sections 4 and 5 that the converse is true for von Neumann algebras and measure algebras. For approximate strong Connes-amenability, the corresponding question is answered in the next theorem which is the approximate version  of [18, Theorem 4.7]. First we need to give a precise definition of this new concept.

\textbf{Definition.} $\A$ is called approximately strongly Connes-amenable if for each unital
Banach $\A$-bimodule $\X$, every $\omega^{*}-\omega^{*}$ continuous
derivation $D\in Z^1(\A ,\X ^*)$ whose range consists of normal elements is approximately inner.

\textbf{Theorem 3.2.} The following conditions are equivalent.

(i) $\A$ has an  approximate normal, virtual diagonal.

(ii) $\A$ is approximately strongly Connes-amenable.

{\it Proof.} $(i)\Longrightarrow(ii)$. This is similar to Theorem 3.1.

$(ii)\Longrightarrow (i)$  Since  $\Delta_{\omega^{*}}$ is $\omega^{*}-\omega^{*}$ continuous then $ker\Delta_{\omega^{*}}$ is $\omega^{*}-$closed and
$$
(L_{\omega^{*}}^{2}(\A,\Bbb C)/^{\perp}ker\Delta_{\omega^{*}})^{*}=ker\Delta_{\omega^{*}}.
$$
So $ker\Delta_{\omega^{*}}$ is a normal dual $\A$-module and $ad_{e\otimes e}$ attains its values in the normal
elements of $ker\Delta_{\omega^{*}}$. By assumption there exists a net
$(N_{\alpha})\subset ker\Delta_{\omega^{*}}$ such that
$$
ad_{e\otimes e}(a)=\lim_{\alpha}a.N_{\alpha}-N_{\alpha}.a\quad (a\in \A).
$$
Let $M_{\alpha}=e\otimes e-N_{\alpha}$. It follows that
$$
a.M_{\alpha}-M_{\alpha}.a\lra0\quad and\quad \Delta_{\omega^{*}}(M_{\alpha})\lra e\quad(a\in \A).
$$
Therefore $(M_\alpha)$ is an approximate normal virtual diagonal for $\A$. $\qed$

We saw that dual Banach algebras with an approximate normal, virtual diagonal are
approximately Connes-amenable, but the converse is likely to be
false in general. We now modify the definition of approximate normal, virtual
diagonal and obtain the desired characterization of approximate Connes-amenability. Let $\A$ be a dual
Banach algebra with predual $\A_{*}$ and let
$\Delta:\A\widehat{\otimes}\A\lra \A$ be the multiplication
map. From [22, Corollary 4.6], we conclude that $\Delta^{*}$ maps $\A_{*}$
into $\sigma WC((\A\widehat{\otimes}\A)^{*})$. Consequently,
$\Delta^{**}$ induces the homomorphism
$$\Delta_{\sigma WC}:\sigma WC((\A\widehat{\otimes}\A)^{*})^{*}\lra \A.$$
With these preparations made, we can now characterize
approximately Connes-amenable, dual Banach algebras through the existence of
certain approximate normal, virtual diagonals. This is indeed an approximate version of [22, Theorem 4.8].

\textbf{Definition.} An approximate $\sigma WC-$virtual diagonal for $\A$ is a net
$(M_{\alpha})$ in$\quad$ \linebreak
$\sigma WC((\A\widehat{\otimes}\A)^{*})^{*}$ such that
$$
a.M_{\alpha}-M_{\alpha}.a\lra0\quad and \quad \Delta_{\sigma WC}(M_{\alpha})\lra e\quad (a\in \A),
$$
the limits being taken in norm.

\textbf{Theorem 3.3.} The following conditions are equivalent.

(i) $\A$ is approximately Connes-amenable.

(ii) There is an approximate $\sigma WC-$virtual diagonal for $\A$.

{\it Proof.}  $(i)\Longrightarrow (ii)$ The map
$$
D:\A\lra \sigma WC((\A\widehat{\otimes}\A)^{*})^{*},\quad a\longmapsto a\otimes e-e\otimes a
$$
is a well defined bounded derivation, since $\A\widehat{\otimes}\A$ can be  embedded canonically into  $\sigma WC((\A\widehat{\otimes}\A)^{*})^{*}$.
Since the dual module $\sigma WC((\A\widehat{\otimes}\A)^{*})^{*}$ is
normal, then it follows that $D$ is $\omega^{*}-\omega^{*}$-continuous. Clearly $D$
attains its values in the $\omega^{*}-$closed submodule $ker
\Delta_{\sigma WC}$ which is a normal dual Banach $\A$-module. So
there is a net $(N_{\alpha})\subset ker\Delta_{\sigma WC}$ such that
$$
Da=\lim_{\alpha}(a.N_{\alpha}-N_{\alpha}.a)\quad (a\in \A).
$$
Letting $M_{\alpha}=e\otimes e-N_{\alpha}$, we see that it is an
approximate $\sigma WC-$virtual diagonal for $\A$.

$(ii)\Longrightarrow (i)$ Let $\X$  be a normal dual Banach $\A$-bimodule.
By Proposition 2.4 we may assume that $\X$  is unital. Let
$D\in Z^1(\A ,\X)$ be a $\omega^{*}-\omega^{*}$-continuous derivation. Define
$$
\theta_{D}:\A\widehat{\otimes}\A\lra  \X ,\quad a\otimes b\longmapsto a.Db.
$$
By [20, Lemma 4.9],  $\theta_{D}^{*}$ maps  the predual $\X _{*}$ into $\sigma WC((\A\widehat{\otimes}\A)^{*})$. Hence $(\theta^{*}|_{\X_{*}})^{*}$ maps $\sigma WC((\A\widehat{\otimes}\A)^{*})^{*}$ into $\X$. Let $(M_{\alpha})\subset \sigma WC((\A\widehat{\otimes}\A)^{*})^{*}$
be an approximate $\sigma WC-$virtual diagonal for $\A$ and let $x_{\alpha}=(\theta^{*}|_{\X_{*}})^{*}(M_{\alpha})$. Observe that $\A\widehat{\otimes}\A$ is $\omega^*$-dense in $ \sigma WC((\A\widehat{\otimes}\A)^{*})^{*}$. So for every $\alpha$ there is a net $(u_\beta^\alpha)$ in $\A\widehat{\otimes}\A$ such that $M_\alpha=\omega^*-\lim_\beta u^\beta_\alpha$. Suppose $c\in\A$ and $t\in\X_*$ are arbitrary. One can easily check that
\begin{eqnarray}
\aligned &x_\alpha .c=\sigma(\X ,\X_*)-\lim_\beta\theta_D(u^\beta_\alpha).c,\quad\text{and}\\
& c. x_\alpha =\sigma(\X ,\X_*)-\lim_\beta c.\theta_D(u^\beta_\alpha)= (\theta^{*}|_{\X_{*}})^{*}(c.M_{\alpha}). \endaligned\label{8}
\end{eqnarray}

On the other hand by [22, Lemma 4.6],  $\Delta^{*}(\A_{*})\subseteq\sigma WC((\A\widehat{\otimes}\A)^{*})$ and hence
\begin{eqnarray}
\Delta_{\sigma WC}(M_\alpha)=\sigma(\A ,\A_*)-lim_\beta\Delta(u_\alpha^\beta). \label{9}
\end{eqnarray}

Suppose $u^\beta_\alpha=\Sigma_k a_k^{\alpha\beta}\otimes b_k^{\alpha\beta}$. Using identities (8) and (9), we obtain
\begin{eqnarray}
\aligned c.x_{\alpha}-x_{\alpha}.c&=(\theta^{*}|_{\X_{*}})^{*}(c.M_{\alpha})-\omega^{*}-\lim_\beta\Sigma_k a_k^{\alpha\beta}.Db_k^{\alpha\beta}.c\\
&=(\theta^{*}|_{\X_{*}})^{*}(c.M_{\alpha})-\lim_\beta\Sigma_k a_k^{\alpha\beta}.D(b_k^{\alpha\beta}c)+\lim_\beta\Sigma_k a_k^{\alpha\beta}b_k^{\alpha\beta}.Dc\\
&=(\theta^{*}|_{\X_{*}})^{*}(c.M_{\alpha}-M_\alpha .c)+\Delta_{\sigma WC}(M_{\alpha}).Dc.\endaligned\label{10}
\end{eqnarray}
By our assumption and (10), we have
$$
Dc=\lim_{\alpha}(c.x_{\alpha}-x_{\alpha}.c)\quad (c\in \A).
$$
This implies that $\A$ is approximately Connes-amenable.$\qed$

\textbf{Remark 3.4.} In the light of proposition 2.3 if we modify the definition of approximate normal virtual diagonal to the following one, then Theorem 3.1 holds also in the case that $\A$ does not have an identity.

``Let $\A$ be a dual Banach algebra (not necessarily unital). A net  $(M_{\alpha})$ in $L_{\omega^{*}}^{2}(\A^\# ,\Bbb C)^{*}$ is called an approximate normal, virtual diagonal for $\A$ if for every $a\in\A^\#$
$$
a.M_{\alpha}-M_{\alpha}.a\lra0\quad and \quad \Delta_{\omega^{*}}(M_{\alpha})\lra e.
$$

\section{Approximate Connes-amenability of von Neumann algebras}

In this section we prove the ``approximate'' analog of the main result of Effros in [5]. First recall some notations from [5]. Let $\A$ be a von Neumann algebra. We call a map $F\in L^{2}(\A ,\Bbb C)$ reduced if there exist states $p,q\in \A _{*}$ and a constant $K$ such that for every $a,b\in \A $,
$$
|F(a,b)|\leq Kp(aa^{*})^{1/2}q(b^{*}b)^{1/2}.
$$
The set $L_{\omega^{*},0}^{2}(\A ,\Bbb C)$  of all such bilinear functionals is an $\A $-
submodule of $L^{2}(\A ,\Bbb C)$ and $L_{\omega^{*},0}^{2}(\A ,\Bbb C)^{*}$ is a normal dual Banach $\A -$bimodule [5, Lammas 2.1 and 2.2]. Also $\A _{*}\otimes \A _{*}\subseteq L_{\omega^{*},0}^{2}(\A ,\Bbb C)$ and $\A \otimes \A $  is identified with an $\A -$submodule of
$L_{\omega^{*},0}^{2}(\A ,\Bbb C)^{*}$. If $\Delta:\A \otimes \A \lra\A$ is the multiplication
map, then $\Delta^{*}$ maps $\A _{*}$ into $L_{\omega^{*},0}^{2}(\A ,\Bbb C)$ and consequently $\Delta^{**}$ drops to an $\A$-bimodule homomorphism
$\Delta_{\omega^{*},0}:L_{\omega^{*},0}^{2}(\A ,\Bbb C)^{*} \lra \A $.

We need the following Lemma in the proof of the next Theorem.

\textbf{Lemma 4.1.} [5, Lemma 2.3] Suppose $\A $ is a finite or properly
infinite von Neumann algebra. Then there is a
$\omega^{*}-\omega^{*}$ continuous linear $\A-$bimodule map
$$
\Phi:L_{\omega^{*},0}^{2}(\A ,\Bbb C)^{*}\lra L_{\omega^{*}}^{2}(\A ,\Bbb C)^{*}
$$
such that $\Delta_{\omega^{*}}\circ\Phi=\Delta_{\omega^{*},0}.$

\textbf{Theorem 4.2.} A von Neumann algebra $\A $ is
approximately Connes-amenable if and only if it has an approximate
normal virtual diagonal.

{\it Proof.} If $\A $ has an approximate normal virtual diagonal
then  by Theorem 3.1, $\A $ is approximately Connes-amenable.

Conversely, suppose $\A $ is approximately Connes-amenable. The dual Banach $\A-$bimodule $L_{\omega^{*},0}^{2}(\A ,\Bbb C)^{*}$ is normal and hence the bounded derivation $D$ defined by
$$
D:\A \lra L_{\omega^{*},0}^{2}(\A ,\Bbb C)^{*}\quad,\quad a\longmapsto a\otimes e_{\A }-e_{\A }\otimes a
$$
is $\omega^{*}-\omega^{*}$ continuous. Since $\Delta_{\omega^{*},0}$ is $\omega^{*}-\omega^{*}$ continuous, then $ker\Delta_{\omega^{*},0}$ is a $\omega^{*}-$closed submodule of $L_{\omega^{*},0}^{2}(\A ,\Bbb C)^{*}$ and we have a Banach $\A-$bimodule isomorphism
$$
(L_{\omega^{*},0}^{2}(\A ,\Bbb C)/^{\perp}ker\Delta_{\omega^{*},0})^{*}\cong ker\Delta_{\omega^{*},0},
$$
 As a result $ker\Delta_{\omega^{*},0}$ is a normal dual Banach $\A-$bimodule and $D(\A )\subseteq\ker\Delta_{\omega^{*},0}.$ Since $\A $ is approximately Connes-amenable, then there exists a net $(N_{\alpha})\subseteq\ker\Delta_{\omega^{*},0}$ such that
$$
Da=\lim_{\alpha}(a.N_{\alpha}-N_{\alpha}.a)\quad(a\in \A ).
$$
If we set $M_{\alpha}=e_{\A }\otimes e_{\A }-N_{\alpha}$, then
$$
\aligned &\lim_{\alpha}(a.M_{\alpha}- M_{\alpha}.a)=0\quad(a\in\A ),\ \text{and}\\
&\lim_{\alpha}\Delta_{\omega^{*},0}(M_{\alpha})=\Delta_{\omega^{*},0}(e_{\A }\otimes e_{\A })=e_{\A }.\endaligned
$$
If $\A $ is finite or properly infinite then by Lemma 4.1, $\widetilde{M_{\alpha}} =\Phi(M_{\alpha})\in
L_{\omega^{*}}^{2}(\A ,\Bbb C)^{*}$ is an approximate normal virtual diagonal because of the following identities
$$\aligned &\lim_{\alpha}(a.\widetilde{M_{\alpha}}-\widetilde{M_{\alpha}}.a)=\lim_{\alpha}
\Phi(a.M_{\alpha}-M_{\alpha}.a)=0\quad(a\in \A ),\\
&\lim_{\alpha}\Delta_{\omega^{*}}(\widetilde{M_{\alpha}})=\lim_{\alpha}\Delta_{\omega^{*}}\circ
\Phi(M_{\alpha})=\lim_{\alpha}\Delta_{\omega^{*},0}(M_{\alpha})=e_{\A }.\endaligned
$$
In the general case, there are central projections $p_{1},p_{2}\in \A $, such that $e_{\A }=p_{1}+p_{2}$,
 $p_{1}\A $ is a finite von Neumann algebra and $p_{2}\A $ is a properly infinite von Neumann
 algebra.  Since $\A $ is approximately Connes-amenable and $\A =p_{1}\A \oplus p_{2}\A$ , then it is easy to see that the von Neumann algebras $p_{1}\A $ and $p_{2}\A $ are approximately Connes-amenable. Therefore
there exist nets $(M_{\alpha})\subseteq L_{\omega^{*}}^{2}(p_{1}\A ,\Bbb C)^{*}$ and
$(M_{\beta})\subseteq L_{\omega^{*}}^{2}(p_{2}\A ,\Bbb C)^{*}$ such that,
$$
\aligned &\lim_{\alpha}(a.M_{\alpha}-M_{\alpha}a)=0\quad(a\in p_{1}\A )\ \text{and}\
\Delta_{\omega^{*}}(M_{\alpha})\lra p_{1},\quad (1)\\
&\lim_{\beta}(a.M_{\beta}-M_{\beta}a)=0\quad(a\in p_{2}\A )\ \text{and}\ \Delta_{\omega^{*}}(M_{\beta})\lra p_{2}.\quad (2)\endaligned
$$
For each $F\in L_{\omega^{*}}^{2}(\A ,\Bbb C)$ define
$$
F_{i}(a,b)=F(a,b)\quad and \quad a,b\in p_{i}\A .\quad(i=1,2)
$$
Clearly $F_{i}\in L_{\omega^{*}}^{2}(p_{i}\A ,\Bbb C)$. Now define the net $(M_{(\alpha,\beta)})\subseteq
L_{\omega^{*}}^{2}(\A ,\Bbb C)^{*}$ by $M_{(\alpha,\beta)}=M_{\alpha}\oplus M_{\beta}$. Then
$$
\langle M_{(\alpha,\beta)},F\rangle=\langle M_{\alpha}\oplus M_{\beta},F\rangle=
\langle M_{\alpha},F_{1}\rangle+\langle M_{\beta},F_{2}\rangle.
$$
For each $c\in \A $ and $a,b\in p_{i}\A ,$
$$
\aligned &(Fc)_{i}(a,b)=(Fc)(a,b)=F(c_{i}a,b)=F_{i}c_{i}(a,b),\\
&(cF)_{i}(a,b)=(cF)(a,b)=F(a,bc_{i})=c_{i}F_{i}(a,b).\endaligned
$$
where $c=c_{1}+c_{2}$, $c_{i}\in p_{i}\A, \ i=1,2.$  Hence
$$
\aligned\langle c.M_{(\alpha,\beta)},F\rangle&=\langle M_{\alpha}\oplus M_{\beta},Fc\rangle\\
&=\langle M_{\alpha},(Fc)_{1}\rangle+\langle M_{\beta},(Fc)_{2}\rangle\\
&=\langle M_{\alpha},F_{1}c_{1}\rangle+\langle M_{\beta},F_{2}c_{2}\rangle\\
&=\langle c_{1}.M_{\alpha},F_{1}\rangle+\langle c_{2}.M_{\beta},F_{2}\rangle\\
&=\langle c_{1}.M_{\alpha}\oplus c_{2}.M_{\beta},F\rangle.\endaligned
$$
Therefore $c.M_{(\alpha,\beta)}=c_{1}.M_{\alpha}\oplus c_{2}.M_{\beta}.$ Similarly
$M_{(\alpha,\beta)}.c=M_{\alpha}.c_{1}\oplus M_{\beta}.c_{2}.$ By (1) and (2) we have,
$$
\aligned\lim_{(\alpha,\beta)}(c.M_{(\alpha,\beta)}-M_{(\alpha,\beta)}.c)&
=\lim_{(\alpha,\beta)}(c_{1}.M_{\alpha}\oplus c_{2}.M_{\beta}-M_{\alpha}.c_{1}\oplus M_{\beta}.c_{2})\\
&=\lim_{(\alpha,\beta)}((c_{1}.M_{\alpha}-M_{\alpha}.c_{1})\oplus(c_{2}.M_{\beta}-M_{\beta}.c_{2}))=0.\endaligned
$$
For $a_{*}\in\A _{*}$ define $a_{*,i}\in (p_{i}\A )_{*}$ by $a_{*,i}(b)=a_{*}(b)$ $(b\in p_{i}\A ).$
So for each $a,b\in p_{i}A,$
$$
\Delta^{*}(a_{*,i})(a,b)=a_{*,i}(ab)=a_{*}(ab)=\Delta^{*}(a_{*})(a,b)=(\Delta^{*}(a_{*}))_{i}(a,b).
$$
Thus,
$$
\aligned\langle \Delta_{\omega^{*}}(M_{(\alpha,\beta)}),a_{*}\rangle&
=\langle M_{\alpha}\oplus M_{\beta},\Delta^{*}(a_{*})\rangle\\
&=\langle M_{\alpha},\Delta^{*}(a_{*})_{1}\rangle+\langle M_{\beta},\Delta^{*}(a_{*})_{2}\rangle\\
&=\langle M_{\alpha},\Delta^{*}(a_{*,1})\rangle+\langle M_{\beta},\Delta^{*}(a_{*,2})\rangle\\
&=\langle \Delta_{\omega^{*}}(M_{\alpha}),a_{*,1}\rangle+\langle \Delta_{\omega^{*}}(M_{\beta}),a_{*,2}\rangle\\
&=\langle \Delta_{\omega^{*}}(M_{\alpha})\oplus \Delta_{\omega^{*}}(M_{\beta}),a_{*,1}\oplus a_{*,2}\rangle.\endaligned
$$
As a result
$\Delta_{\omega^{*}}(M_{(\alpha,\beta)})=\Delta_{\omega^{*}}(M_{\alpha})\oplus
\Delta_{\omega^{*}} (M_{\beta})$. Also
$$
\Delta_{\omega^{*}}(M_{(\alpha,\beta)})-e_{\A }=(\Delta_{\omega^{*}}(M_{\alpha})-p_{1})
\oplus (\Delta_{\omega^{*}}(M_{\beta})-p_{2}).\quad (3)
$$
Finally based on (1),(2) and (3) we have
$$
\lim_{(\alpha,\beta)}\Delta_{\omega^{*}}(M_{(\alpha,\beta)})\lra e_{\A }.
$$
It follows that net $(M_{(\alpha,\beta)})$ is an approximate normal virtual diagonal for $\A $. $\qed$

\section{Approximate Connes-amenability of $M(G)$}

 In this section we characterize approximate Connes-amenable measure algebras on locally compact groups. Throughout this section $G$ is a locally compact group, $G^{op}$ denotes the same group, with reversed multiplication and $\Delta:M(G)\widehat{\otimes}M(G)\lra M(G)$
 is the multiplication map. We recall some terminology from [20]. A bounded function $f:G\times G^{op}\lra \Bbb C$ is called separately $C_{0}$ if for each $x\in G$, the function
  $G^{op}\lra \Bbb C$, $y\mapsto f(x,y)$, belongs to $C_{0}(G^{op})$, and for each
$y\in G^{op}$, the function $G\lra \Bbb C$, $x\mapsto f(x,y)$ belongs to $C_{0}(G)$. The collection of all separately $C_{0}-$functions is denoted by  $SC_{0}(G\times G^{op})$. Let $LUC(G)$ be the commutative $C^*$-algebra of left uniformly continuous functions on $G$ and $G_{LUC}$ be its character space.
The set
$$
\{ f\in LUC(G\times G^{op})\ :\ \phi .f\in SC_0(G\times G^{op})\ \text{for} \ \text{all}\ \phi\in (G\times G^{op})_{LUC}\}
$$
which is denoted by $LUCSC_{0}(G\times G^{op})$ is a closed $M(G)$-submodule of $SC_{0}(G\times G^{op})$ whose dual  is a normal dual Banach $M(G)$-bimodule.
Moreover $\Delta_{*}=\Delta^{*}|_{C_{0}(G)}$ maps $C_{0}(G)$ into $LUCSC_{0}(G\times G^{op})$ [20, theorem 4.4]. Therefore $\Delta^{**}$ turns into an $M(G)$-bimodule homomorphism $\widetilde{\Delta}:LUCSC_{0}(G\times G^{op})^{*}\lra M(G).$

\textbf{Proposition 5.1.} If $M(G)$ is approximately Connes-amenable, then $G$ is amenable.

{\it Proof}.  First we show that here is a net $(M_{\alpha})\subseteq LUCSC_{0}(G\times G^{op})^{*}$ such that
$$
\mu.M_{\alpha}-M_{\alpha}.\mu\lra 0\quad (\mu\in M(G))\quad \text{and}\quad
\widetilde{\Delta}(M_{\alpha})\lra \delta_{e}.
$$
It is easy to see that the map
$$
D:M(G)\lra LUCSC_{0}(G\times G^{op})^{*},\quad \mu\longmapsto \mu\otimes \delta_{e}-\delta_{e}\otimes\mu
$$
is a bounded derivation. By [20, Proposition 3.2] the maps $\mu\longmapsto \mu\otimes \delta_{e}$ and \linebreak
$\mu\longmapsto\delta_{e}\otimes \mu$ form $M(G)$ into $M(G\times G^{op})$ are $\omega^{*}-\omega^{*}$continuous and hence so is $D$. Moreover $D(M(G))\subseteq ker \widetilde{\Delta}.$ Since
 $\widetilde{\Delta}$ is a $\omega^{*}-\omega^{*}$ continuous bimodule homomorphism, then
 $ker \widetilde{\Delta}$ is a $\omega^{*}-$ closed submodule and
$$(LUCSC_{0}(G\times G^{op})/^{\perp}ker\widetilde{\Delta})^{*}\cong ker\widetilde{\Delta}$$
as Banach $M(G)$-bimodules. By [20, Theorem 4.4(ii)], $LUCSC_{0}(G\times G^{op})^{*}$ is a
normal dual Banach $M(G)$-bimodule, and so is $ker\widetilde{\Delta}.$ Since $M(G)$ is
approximately Connes-amenable, then there is a net $(N_{\alpha})\subseteq ker\widetilde{\Delta}$
such that $D\mu=\lim_{\alpha}\mu.N_{\alpha}-N_{\alpha}.\mu\quad(\mu\in M(G))$. The identities
$\mu.N_{\alpha}=(\mu\otimes \delta_{e})*N_{\alpha}$ and $N_{\alpha}.\mu=N_{\alpha}*(\delta_{e}\otimes \mu)$ imply that if we set $M_{\alpha}=\delta_{e}\otimes \delta_{e}-N_{\alpha}$, then the net
 $(M_{\alpha})$ has the required properties.

Since $\widetilde{\Delta}(M_{\alpha})\lra \delta_{e}$, then we can suppose that $M_{\alpha}\neq 0$ for every $\alpha$ and if we consider $M_{\alpha}$ as a measure on the character space of the
commutative $C^{*}$-algebra $LUCSC_{0}(G\times G^{op})$, then the total variation
$|M_{\alpha}|$ is a non-zero element of $LUCSC_{0}(G\times G^{op})^{*}$. Observe that
$$
|\delta_{g}.M_{\alpha}|=|(\delta_{g}\otimes \delta_{e})*M_{\alpha}|=(\delta_{g}\otimes \delta_{e})*|M_{\alpha}|=\delta_{g}.|M_{\alpha}|,
$$
and similarly $|M_{\alpha}.\delta_{g}|=|M_{\alpha}|.\delta_{g}.$ Thus
$$
\|\delta_{g}.|M_{\alpha}|-|M_{\alpha}|.\delta_{g}\|=\||\delta_{g}.M_{\alpha}|-|M_{\alpha}.\delta_{g}|\|\leq\|M_{\alpha}.\delta_{g}-\delta_{g}.M_{\alpha}\|\lra 0.
$$
On the other hand the convergence $\widetilde{\Delta}(M_{\alpha})\lra \delta_{e}$ implies that the net $(1/\|M_{\alpha}\|)$ is bounded and so if we set $N_{\alpha}=|M_{\alpha}|/\|M_{\alpha}\|$, then $\delta_{g}.N_{\alpha}-N_{\alpha}.\delta_{g}\lra 0.$ Let $N$ be a $\omega^{*}-$ cluster point of  $(N_{\alpha})$. Then
$\delta_{g}.N=N.\delta_{g}$ for every $g\in G$. As in the proof of [20, Theorem 5.3], $LUC(G\times G^{op})$ can be considered as a $C^{*}-$subalgebra of $LUCSC_{0}(G\times G^{op})^{**}$; So  in particular
$\langle f,N\rangle$ is well defined for each $f\in LUC(G\times G^{op}).$ Note that  the embedding
 of $LUC(G\times G^{op})$ into $LUCSC_{0}(G\times G^{op})^{**}$ is an $M(G)-$bimodule
 homomorphism. Define
$$
m:LUC(G)\lra C\quad, \quad f\longmapsto\langle N,f\otimes1\rangle.
$$
Since $f\otimes1\in LUC(G\times G^{op})$, then $m$ is a well-defined, positive, linear functional
whose normalization is a left invariant mean on $LUC(G)$  as in the proof of [20, Theorem 5.3]. Therefore $G$ is amenable. $\qed$

Combination of the preceding proposition, Theorem 3.1 and the main result of [21] leads to the following result.

\textbf{Theorem 5.2.} The following conditions are equivalent.

(i) $G$ is amenable.

(i) $M(G)$ is approximately Connes-amenable.

(iii) $M(G)$ is  Connes-amenable.

(iv) $M(G)$ has an approximate normal, virtual diagonal.

The algebra $WAP(G)$ of weakly almost periodic functions on $G$ is a commutative $C^{*}-$algebra which is a left introverted subspace of $L^{\infty}(G)$ [24, Lemma 6.3]. Thus $WAP(G)^*$ is a dual Banach algebra which is identified with $WAP(L^\infty(G))^*$. In the next proposition we identify the relationship between approximate amenability of $\A$ and approximate Connes-amenability of $WAP(\A^{*})^{*}$ in the special case of group algebras.

\textbf{Proposition 5.3.} $G$ is  amenable if and only if $WAP(L^\infty(G))^{*}$ is approximately Connes-amenable.

{\it Proof.} Suppose $G$ is amenable. By [14, Theorem 2.5] and [8, Theorem
3.2], $G$ is amenable if and only if $L^{1}(G)$ is approximately
amenable. Since the image of $L^{1}(G)$ is $\omega^{*}-$dense in
WAP$(L^\infty(G))^{*}$, then WAP$(L^{1}(G)^{*})^{*}$ is Connes-amenable [18, Proposition 4.2(i)].

Conversely suppose $WAP(L^\infty(G))^{*}$ is approximately Connes-amenable. Since
 \linebreak
 $C_{0}(G)\subseteq WAP(G)$, the restriction map from WAP$(G)^{*}$
onto $M(G)$ is a $\omega^{*}-\omega^{*}$ continuous algebra homomorphism. Consequently $M(G)$ is approximately Connes-amenable and by Theorem 5.2, $G$ is amenable. $\qed$

\section{APPROXIMATE CONNES AMENABILTIY OF $\A^{**}$}

If $\A$ is a dual Banach algebra such that $\A^{**}$ is Connes-amenable, then so is $\A$ [4]. In the following theorem we extend this result to approximate Connes-amenability.

\textbf{Theorem 6.1.} Let $\A$ be an Arens regular Banach algebra such that $\A^{**}$ is approximately Connes-amenable.

(i) If $\A$ is a dual Banach algebra, then $\A$ is approximately Connes-amenable.

(ii) If $\A$ is an ideal in $\A^{**}$ and $\A^{**}$  has an identity then $\A$ is approximately amenable.

{\it Proof.} (i) Suppose $\X$  is a normal dual Banach $\A$-bimodule, and
$\pi:\A^{**}\lra \A$ is the restriction map to $\A_*$. Then $\pi$ is a
$\omega^{*}-\omega^{*}$ continuous homomorphism. Therefore $\X$  is a
normal dual Banach $\A^{**}$-bimodule with the following actions
$$a^{**}.x=\pi(a^{**})x\quad ,\quad x.a^{**}=x\pi(a^{**})\quad (x\in  \X ,a^{**}\in \A^{**}).$$

Let $D:\A\lra  \X $ be a $\omega^{*}-\omega^{*}$ continuous
derivation. It is easy to see that $Do\pi:\A^{**}\lra  \X $
is a $\omega^{*}-\omega^{*}$ continuous derivation. Since $\A^{**}$
is approximately Connes-amenable, than there exists a net
$(x_{\alpha})\subseteq  \X $ such that
$$Do\pi(a^{**})=\lim_{\alpha}a^{**}.x_{\alpha}-x_{\alpha}.a^{**}\quad (a^{**}\in \A^{**}).$$
So
$$D(a)=\lim_{\alpha}a.x_{\alpha}-x_{\alpha}.a\quad (a\in \A).$$

(ii) By [8, Proposition 2.5] in order to show that $\A$ is approximately
amenable it is sufficient to show that every $D\in Z^{1}(\A, \X ^{*})$ is
approximately inner for each neo-unital Banach $\A$-module.

Let $\X$  be a neo-unital Banach $\A$-bimodule, and let $D\in
Z^{1}(\A, \X ^{*})$. As in the proof of [19, Theorem 4.4.8] one can show that $\X ^{*}$ is a normal dual Banach
$\A^{**}$-bimodule and D has a unique extension $\widetilde{D}\in
Z^{1}(\A^{**}, \X ^{*})$. From the approximate Connes-amenability of
$\A^{**}$ we conclude that $\widetilde{D}$, and hence D is inner.
It follows that $\A$ is approximately amenable.$\qed$

\textbf{Theorem 6.2.} Suppose $\A$ is a Banach algebra with a bounded approximate identity $(e_\beta)$ and $B(\A ,\A^*)=W(\A ,\A^*)$. If $\A^{**}$ is approximately strongly Connes-amenable, then $\A$ is approximately amenable.

{\it Proof.} Following the argument of [18, Theorem 4.8] we see that $\A$ is Arens regular and hence $\A^{**}$ is a dual Banach algebra. Moreover
$$
(\A\widehat{\otimes} \A)^{**}\cong L^2_{\omega^*}(\A^{**},\Bbb C)^*
$$
as Banach $\A$-bimodules. Since $\A$ is Arens regular and has a bounded approximate identity, then $\A^{**}$ has an identity $e$. Thus by Theorem 3.2 $\A^{**}$ has  an approximate normal virtual diagonal $(M_\alpha)\subset L^2_{\omega^*}(\A^{**},\Bbb C)^*$. Now set $M^{''}_{(\alpha, \beta)}=M_\alpha+e_\beta\otimes e_\beta$
and $F_{(\alpha, \beta)}=G_{(\alpha, \beta)}=e_\beta$.  Then for every $a\in\A$ we have
$$
\aligned aM^{''}_{(\alpha, \beta)}-M^{''}_{(\alpha, \beta)}a&+F_{(\alpha, \beta)}\otimes a-a\otimes G_{(\alpha, \beta)}\\
&=aM_\alpha-M_\alpha a+ae_\beta\otimes e_\beta-e_\beta\otimes e_\beta a+e_\beta\otimes a-a\otimes e_\beta\\
&=(aM_\alpha-M_\alpha a)+(ae_\beta-a)\otimes e_\beta+e_\beta\otimes (a-e_\beta a)\lra 0.\endaligned
$$
Moreover $aF_{(\alpha, \beta)}\lra a,\quad G_{(\alpha, \beta)}a\lra a$ and
$$
\Delta^{**}(M^{''}_{(\alpha, \beta)})a-F_{(\alpha, \beta)}a-G_{(\alpha, \beta)}a=\Delta^{**}(M_\alpha)a+e^2_\beta a-e_\beta a-e_\beta a\lra 0.
$$
Therefore by [8, Corollary 2.2] $\A$ is approximately amenable. $\qed$

{\bf Acknowledgments.} The authors would like to express their sincere thanks to Professor F. Ghahramani and Professor A. T. M. Lau for their valuable comments.


\begin{thebibliography}{99}

\bibitem{ref1}%1
H. G. Dales, Banach algebras and automatic continuity, Clarendon Press, Oxford, 2000.

\bibitem{ref1}%1
H. G. Dales, F. Gahramani, and N. Gronbeak, Derivations into
iterated duals of Banach algebras,  Studia Math. 128 (1)
(1998), 19-53.

\bibitem{ref1}%1
H. G. Dales, A. T. M. Lau and D. Strauss, Banach algebras on semigroups and their compactifications,  Memoir Amer. Math. Soc. 966 (2010).

\bibitem{ref1}%1
M. Daws, Connes-amenability of bidual and weighted semigroup algebras, Math. Scand. 99 (2006), no.2, 217-246.


\bibitem{ref1}%1
E. G. Effros, Amenability and virtual diagonals for von Neumann
algebras, J. Funct. Anal. 78 (1988), 137-153.

\bibitem{ref1}%1
G. H. Esslamzadeh and B. Shojaee, Approximate weak amenability of Banach algebras, To appear.

\bibitem{ref1}%1
F. Ghahramani and A. T. M. Lau, Approximate weak amenability, derivations and Arens regularity of Segal algebras, Studia Math. 169 (2005), 189-205.

\bibitem{ref1}%1
F. Ghahramani and R. J. Loy, Generalized notions of amenability, J. Funct. Anal. 208 (2004),229-260.

\bibitem{ref1}%1
F. Ghahramani, R. J. Loy and Y. Zhang, Generalized notions of amenability II, J. Funct. Anal. 254 (2008), 1776-1810.

\bibitem{ref1}%1
F. Gourdeau, Amenability and the second dual of Banach algebras, Studia Math.  125 (1997), 75-81.

\bibitem{ref1}%1
U. Haagerup, All nuclear $C^*$-algebras are amenable, Invent. Math. 74 (1983), 305-319.

\bibitem{ref1}%1
A. Ya. Helemskii, The Homology of Banach and Topological Algebra (translated from the Russian). Kluwer Academic Publishers, 1989.

\bibitem{ref1}%1
A. Ya. Helemskii, Homlogical essence of amenability in the sence of A. Connes: the injectivity of the predual bimodule (translated from the Russion). Math. USSR-Sb 68 (1991), 555-566.


\bibitem{ref1}%1
B. E. Johnson, Cohomology in Banach algebras, Mem. Amer. Math. Soc. 127 (1972).

\bibitem{ref1}%1
B. E. Johnson, Approximate diagonals and cohomology of certain annihilator Banach algebras, Amer. J. Math. 94 (1972), 685-698.

\bibitem{ref1}%1
B. E. Johnson, R. V. Kadison, and J. Ringrose, Cohomology of operator algebras, III. Bull. Soc. Math. France 100 (1972), 73-79.

\bibitem{ref1}%1
A. T. M. Lau and R. Loy, Weak amenability of Banach algebras on locally compact groups, J. Funct. Anal. 145 (1997), 175-204.

\bibitem{ref1}%1
V. Runde, Amenability for dual Banach algebras, Studia Math. 148 (2001), 47-66.

\bibitem{ref1}%1
V. Runde, Lectures on amenability, Lecture Notes in Mathematics 1774, Springer-Verlage, Berlin, 2002.

\bibitem{ref1}%1
V. Runde, Connes-amenability and normal virtual diagonals for
measure alebras I, J. London Math. Soc. 67(2003), 643-656.

\bibitem{ref1}%1
V. Runde, Connes-amenability and normal virtual diagonals for
measure algebras II, Bull Aust. Math. Soc. 68 (2003), 325-328.

\bibitem{ref1}%1
V. Runde, Dual Banach algebras: Connes-amenability, normal, virtual
diagonal, and injectivity of the predual bimodule, Math. Scand. 95 (2004), 124-144.


\bibitem{ref1}%1
B. Shojaee, G. H. Esslamzadeh and A. Pourabbas, First order cohomology of $\ell^1$-Munn algebras, Bull. Iranian Math. Soc. 35 (2009), 211-219.

\bibitem{ref1}%1
J.C.S. Wong, Topologically stationary locally compact groups and amenability, Trans. Amer. Math. Soc. 144(1969), 351-363.

\end{thebibliography}
\end{document}